\documentclass{amsart}

\usepackage{enumerate}

\newtheorem{theorem}{Theorem}
\newtheorem{proposition}[theorem]{Proposition}
\newtheorem{lemma}[theorem]{Lemma}

\numberwithin{equation}{section}

\newcommand{\Ot}[1]{\mathcal O(t^{#1})}
\newcommand{\Os}[1]{\mathcal O(s^{#1})}
\newcommand{\Odt}[1]{\mathcal O(\Delta t)^{#1}}

\newcommand{\pd}[2]{\frac{\partial {#1}}{\partial {#2}}}
\newcommand{\pdd}[3]{\frac{\partial^2 {#1}}{\partial {#2} \partial {#3}}}

\begin{document}

\title[Unexpected local extrema]
      {Unexpected local extrema \\ for the Sendov conjecture}
\author{Michael J. Miller}
\address{Department of Mathematics\\Le Moyne College\\Syracuse, New 
York 13214\\USA}
\email{millermj@lemoyne.edu}
\thanks{Last revised 9-Aug-2007} 

\begin{abstract}
Let $S(n)$ be the set of all polynomials of degree $n$ with all roots
in the unit disk, and define $d(P)$ to be the maximum of the distances
from each of the roots of a polynomial $P$ to that root's nearest
critical point.  In this notation, Sendov's conjecture asserts that
$d(P) \le 1$ for every $P \in S(n)$.

Define $P \in S(n)$ to be \emph{locally extremal} if $d(P) \ge d(Q)$
for all nearby $Q \in\nobreak S(n)$, and note that maximizing $d(P)$ over all
locally extremal polynomials~$P$ would settle the Sendov conjecture.

Prior to now, the only polynomials known to be locally extremal were
of the form $P(z)=c(z^n+a)$ with $|a|=1$.
In this paper, we determine sufficient conditions for real polynomials
of the form 
\begin{displaymath}
P(z)=\int_\beta^z (w-a)^{n-3} (w^2+bw+c)\, dw \text{\quad with $0<\beta<1$}
\end{displaymath}
to be locally extremal, and we use these conditions to find locally
extremal polynomials of this form of degrees $n=8$, $9$, $12$, $13$,
$14$, $15$, $19$, $20$, and~$26$.
\end{abstract}

\maketitle

\section{Introduction}\label{section1} 

In 1958, Sendov conjectured that if a polynomial (with complex
coefficients) has all its roots in the unit disk, then within one unit
of each of its roots lies a root of its derivative.  This conjecture
has yet to be settled, although it has been the subject of more than
80 papers over the intervening years, and so has been verified for
many special cases.  These have been documented by Sendov
\cite{Sendov}, Schmeisser~\cite{Schmeisser}, Sheil-Small
\cite[Chapter~6]{Sheil-Small} and Rahman and Schmeisser
\cite[Section~7.3]{Rahman-Schmeisser}.

Let $n\ge2$ be an integer and let $S(n)$ be the set of all complex
polynomials of degree $n$ with all roots in the unit disk.
For a polynomial $P$ with roots $z_1, \dots, z_n$ and critical points
$\zeta_1, \dots, \zeta_{n-1}$, define
\begin{displaymath}
   d(P) = \max_{1 \le i \le n} 
      \left\{ \min_{1 \le j \le n-1} |z_i-\zeta_j| \right\}.
\end{displaymath}
If $P \in S(n)$, then the Gauss-Lucas Theorem \cite[Theorem
2.1.1]{Rahman-Schmeisser} implies that each $d(P) \le 2$, and Sendov's
conjecture asserts that each $d(P) \le 1$.

We will say that a polynomial $P$ is \emph{expected} if it is of the
form $P(z)=c(z^n+a)$ with $|a|=1$.  In 1972, Phelps and Rodriguez
defined a polynomial $P \in S(n)$ to be \emph{extremal} if $d(P)=\sup
\{ d(Q) : Q \in S(n) \}$, and conjectured \cite[after
Theorem~5]{Phelps-Rodriguez} that extremal polynomials are all
expected.  Since any expected polynomial $P$ has $d(P)=1$, this
conjecture implies Sendov's conjecture.  The Phelps-Rodriguez
conjecture has also been verified for a number of special cases, as
documented by Rahman and Schmeisser \cite[Section
7.3]{Rahman-Schmeisser}.

Define an \emph{$\epsilon$-neighborhood} of $P \in S(n)$ to be the set
of all polynomials $Q \in~S(n)$ whose roots are within $\epsilon$ of
the roots of $P$ (in the sense that the roots of~$Q$ can be paired
with the roots of $P$ so that in each pair, the difference of the
roots has a modulus less than $\epsilon$).  Define a polynomial $P \in
S(n)$ to be \emph{locally extremal} if $d(P) \ge d(Q)$ for all $Q$ in
some $\epsilon$-neighborhood of $P$, and note that maximizing
$d(P)$ over all locally extremal polynomials $P$ would settle the
Sendov conjecture.

The expected polynomials are all locally extremal, as was demonstrated
by V\^aj\^aitu and Zaharescu~\cite{Vajaitu-Zaharescu} and by Miller
\cite[Theorem~3]{Miller}.  Given this, it is tempting to approach
Sendov's conjecture by trying to show that all locally extremal
polynomials must be expected.  Indeed, Schmieder has made several such
attempts~\cite{Schmieder}, although Borcea has shown
\cite[section~1]{Borcea-2} that each contains a technical flaw.

In this paper, we prove
\begin{theorem}\label{theorem1}
For each $n \in \{8, 9, 12, 13, 14, 15, 19, 20, 26 \}$, there are
locally extremal polynomials $P \in S(n)$ of the form
\begin{displaymath}
P(z)=\int_\beta^z (w-a)^{n-3} (w^2+bw+c)\, dw \text{\quad with $0<\beta<1$}.
\end{displaymath}
\end{theorem}
\noindent Note that polynomials of this form would not be expected, as
each has a root~$\beta$ that is not on the unit circle.  Thus
Theorem~\ref{theorem1} implies that there are locally extremal
polynomials that are not expected.

Theorem \ref{theorem1} has several important consequences.  First, it
shows that the technical flaws in Schmieder's approach cannot be
patched.  Second, the structure of these polynomials (in particular,
the multiple critical point) identifies circumstances that potential
variational proofs of Sendov's conjecture will need to address.
Finally, finding these polynomials is a significant new step in
identifying all local extrema.

In section \ref{section2}, we list eight properties that are sufficient
for a polynomial to be locally extremal, and in sections \ref{section3}
and \ref{section4} we verify that these properties suffice.  In
section \ref{section5}, we describe how to construct polynomials that
satisfy all these properties, and in section \ref{section6} we list
the resulting polynomials, thereby verifying Theorem~\ref{theorem1}.

\section{Properties}\label{section2} 

In section 4 we prove that for a real polynomial $P$ of degree $n \ge
5$ to be locally extremal, it suffices for it to satisfy 8 properties,
beginning with the following.
\begin{enumerate}[\bf A:]
   \item \label{A}All roots of $P$ lie in the closed unit disk.
   \item \label{B}All roots of $P$ that are on the unit circle are simple.
   \item \label{C}$P$ has a root at $\beta$, with $0 < \beta < 1$.
   \item \label{D}All critical points of $P$ lie on a circle of positive
                radius centered at $\beta$.
   \item \label{E}$P$ has a real critical point $a<\beta$ of order $n-3$.
\end{enumerate}

Let $z_1, \dots, z_n$ be the roots of $P$, numbered so that
$z_1,\dots,z_m$ are on the unit circle.  Let $\zeta_1, \dots,
\zeta_{n-1}$ be the critical points of $P$, numbered (as allowed by
property~\ref{E}) so that $\zeta_j=a$ for $j \ge 3$.

Note that property~\ref{D} implies that $P$ has a simple root at $\beta$.
Our next property is
\begin{enumerate}[\bf A:]
   \setcounter{enumi}{5}
   \item \label{F} For each $z_i \ne \beta$, we have 
                       $\min_{1 \le j \le n-1} |z_i-\zeta_j| 
                       < \min_{1 \le j \le n-1} |\beta-\zeta_j| < 1$.         
\end{enumerate}

To examine the effects of changing $\beta$ and the $\zeta_j$ by small
amounts, denoted by $\Delta\beta$ (which we will require to be real)
and $\Delta\zeta_j$, we will use the following notation.
\begin{displaymath}\label{2.1}\tag{2.1}
\begin{aligned}
E_k &= -\Re \left[ \frac{\Delta\zeta_k-\Delta\beta}{\zeta_k-\beta} \right]
          \text{\quad for $k=1$ and $k=2$,} \\
F_k &= 0  \text{\quad for $k=1$ and $k=2$,} \\
E_3 &= \frac{\Re[\sum_{j=3}^{n-1} \Delta\zeta_j] - (n-3) \Delta\beta}
           {\beta-a}, \text{ and} \\
F_3 &= - \frac{\sum_{j=3}^{n-1} (\Im[\Delta\zeta_j])^2}{2(\beta-a)^2}. \\
\end{aligned}
\end{displaymath}

Since (up to a constant multiple) $P'(z)= \prod_{j=1}^{n-1}
(z-\zeta_j)$ and by property~\ref{C} we have $P(z)=\int_\beta^z
P'(w)\,dw$, then the roots of $P$ are functions of $\beta$ and the
$\zeta_j$.  By property~\ref{B}, the roots of $P$ that are on the unit
circle (being simple) are differentiable functions of $\beta$ and the
$\zeta_j$.

Recall that $z_1,\dots,z_m$ are on the unit circle.  For $i=1, \dots,
m$ define
\begin{align*}
E_{i+3} &= \Re \left[ \frac{1}{z_i} \pd{z_i}{\beta} \right] \Delta\beta
   + \sum_{j=1}^2 \Re \left[ \frac{1}{z_i} \pd{z_i}{\zeta_j}\Delta\zeta_j
                      \right]
  + \Re \left[ \frac{1}{z_i} \pd{z_i}{\zeta_3} \sum_{j=3}^{n-1} \Delta\zeta_j
        \right]  \text{ and} \\
F_{i+3} &= - \Re \left[ \frac{1}{2 z_i P'(z_i)} 
                \int_\beta^{z_i} \frac{P'(w)\,dw}{(w-a)^2}
         \right] \sum_{j=3}^{n-1} (\Im[\Delta\zeta_j])^2. 
\end{align*}

Note that for a fixed polynomial $P$, the $m+3$ expressions $E_k$ are
all linear in the 7 real ``variables'' $\Delta\beta$,
$\Re[\Delta\zeta_1]$, $\Im[\Delta\zeta_1]$, $\Re[\Delta\zeta_2]$,
$\Im[\Delta\zeta_2]$, $\Re[\sum_{j=3}^{n-1} \Delta\zeta_j]$ and
$\Im[\sum_{j=3}^{n-1} \Delta\zeta_j]$, and that the $m+3$ expressions
$F_k$ are all constant multiples of the real ``variable''
$\sum_{j=3}^{n-1} (\Im[\Delta\zeta_j])^2$.  Our final two properties
are
\begin{enumerate}[\bf A:]
   \setcounter{enumi}{6}
   \item \label{G}There are constants $c_k>0$ (depending on $P$, but
         independent of the 8 ``variables'') so that the sums 
         $\sum_{k=1}^{m+3} c_k E_k =0$ and $\sum_{k=1}^{m+3} c_k F_k =  
         \sum_{j=3}^{n-1} (\Im[\Delta\zeta_j])^2$.
   \item \label{H}The coefficient matrix of the system 
         $\{E_k =0: k=1,\dots,m+3\}$ in our 7 ``variables'' is of 
         rank 7.
\end{enumerate}

\section{Preliminary calculations} \label{section3} 

We will show that a real polynomial $P$ of degree $n \ge 5$ satisfying
properties \ref{A}--\ref{H} is locally extremal, as follows.

Define $r$ to be the radius of the circle in property~\ref{D} and
recall that (up to a constant multiple) $P(z)=\int_\beta^z \prod_{j=1}^{n-1}
(w-\zeta_j)\,dw$. From our properties \ref{A} and \ref{D} we
know that $P$ has all its roots in the closed unit disk and all its
critical points on a circle of radius $r>0$ centered at $\beta$.
Define the $n$-tuple $(\Delta\beta, \Delta\zeta_1, \dots,
\Delta\zeta_{n-1}) \in C^n$ to be an \emph{improvement} of $P$ if
$\Delta\beta$ is real and if the polynomial
\begin{equation}\label{3.1}
   \int_{\beta+\Delta\beta}^z \prod_{j=1}^{n-1} [w-(\zeta_j+\Delta\zeta_j)]\,dw
\end{equation}
has all its roots in the closed unit disk and all its critical points
strictly outside the circle of radius $r$ centered at
$\beta+\Delta\beta$.  By property \ref{F} we know that $r<1$, so there
is at least one improvement, namely $(1-\beta, -\zeta_1, \dots,
-\zeta_{n-1})$.

For an improvement $I=(\Delta\beta, \Delta\zeta_1, \dots, \Delta\zeta_{n-1})$,
define 
\begin{displaymath}
  ||I|| = \left( |\Delta\beta|^2 
           + \sum_{j=1}^{n-1} |\Delta\zeta_j|^2 \right)^{1/2}.
\end{displaymath}
Note that $||I||>0$, for if $||I||=0$, then the critical points of
\eqref{3.1} would be on (and thus not strictly outside) the circle of
radius $r$ centered at $\beta+\Delta\beta$.

Take any $\epsilon$ with $0<\epsilon<\beta/2$ and recall our
definition of $\epsilon$-neighborhood from section~\ref{section1}.
The critical points of $P$ are continuous functions of the roots of
$P$, so there is a $\delta>0$ such that for all polynomials $Q$ in a
$\delta$-neighborhood of $P$, the roots and critical points of $Q$ are
within $\epsilon$ of the roots and critical points of~$P$.  By
Property~\ref{F}, we know that $d(P)=r$, so if $P$ is not locally
extremal, then there is a polynomial $\hat P \in S(n)$ in the
$\delta$-neighborhood of $P$ with $d(\hat P)>r$.  Let $\hat
\beta$ be the root of $\hat P$ paired with $\beta$ (so
$|\hat\beta-\beta|<\epsilon$), and for $j=1,\dots,n-1$, let
$\hat\zeta_j$ be the critical point of $\hat P$ paired with~$\zeta_j$
(so $|\hat\zeta_j-\zeta|<\epsilon$).  Note that $|\hat\beta| \ge
|\beta| -|\hat\beta-\beta| \ge \beta/2 > 0$, define
$u=|\hat\beta|/\hat\beta$ and note that $|u-1| =
||\hat\beta|-\hat\beta|/|\hat\beta| \le 2\epsilon/(\beta/2)$.  Since
$|u|=1$, then the transformation $z \to uz$ is a rotation about the
origin, so $I=(u\hat\beta-\beta, u\hat\zeta_1-\zeta_1, \dots,
u\hat\zeta_{n-1}-\zeta_{n-1})$ is an improvement of $P$.  Note that
each $|u\hat\zeta_j-\zeta_j| \le |u||\hat\zeta_j-\zeta_j| +
|\zeta_j|(u-1) \le C\epsilon$ for some constant $C$ (and similarly for
$|u\hat\beta-\beta|$), so $||I||$ can be made arbitrarily small.

We have just seen that if $P$ is not locally extremal then there are
improvements $I$ of $P$ with $||I||$ arbitrarily small.  This means
that we can show that $P$ is locally extremal by proving
\begin{theorem}\label{theorem2}
If $P$ is a real polynomial of degree $n \ge 5$ that satisfies
properties \ref{A}--\ref{H}, then there is a constant $A>0$ so that
every improvement $I$ of $P$ has $||I|| \ge A$.
\end{theorem}

We begin with some preliminary calculations.  For $t>0$, we'll say
that a (real or complex) quantity is $\Ot{}$ if its modulus is bounded
by a constant multiple of $t$.  Given this, we have
\begin{lemma}\label{lemma3}
If $z=\Ot{}$ then 
\begin{align*}
   &\text{\rm (1)\qquad}  
      |1+z| = 1 + \Re[z] + \Ot{2}, \text{ and}\\
   &\text{\rm (2)\qquad}
      |1+z| = 1 + \Re[z] + (1/2)(\Im[z])^2 + \Ot{3}.
\end{align*}
\end{lemma}
\begin{proof}
Note that $\Re[z]=\Ot{}$ and $|z|^2=\Ot{2}$.  The results follow from
the equality $|1+z|=(1 + 2\Re[z] + |z|^2)^{1/2}$, the Taylor series
$(1+s)^{1/2} = 1 + s/2 -s^2/8 + \Os{3}$ and the substitution
$s=2\Re[z]+|z|^2$.
\end{proof}

We will be examining the relationships between the roots and critical
points of~$P$.  We calculate the partial derivatives of these
relationships with
\begin{lemma}\label{lemma4}
Write $P'(z)=\prod_{j=1}^{n-1}(z-\zeta_j)$ and suppose that $z_i$ is a simple
root of $P(z)=\int_\beta^z P'(w)\,dw$.  Then
\begin{align*}
   &\text{\rm (1)\quad}  
    \text{$z_i$ is an analytic function of 
         $\beta, \zeta_1, \dots, \zeta_{n-1}$},\\
   &\text{\rm (2)\quad}  
    \pd{z_i}{\beta} =\frac{P'(\beta)}{P'(z_i)},\\
   &\text{\rm (3)\quad}  
    \pd{z_i}{\zeta_j} = 
         \frac{1}{P'(z_i)} \int_\beta^{z_i} \frac {P'(w)\,dw}{w-\zeta_j},\\
   &\text{\rm (4)\quad}  
    \pd{^2 z_i}{\zeta_j^2} 
         = \frac{2}{z_i-\zeta_j} \pd{z_i}{\zeta_j}
         -  \frac{P''(z_i)}{P'(z_i)} \left( \pd{z_i}{\zeta_j} \right)^2 
         \text{\quad and}\\
   &\text{\rm (5)\quad}  
    \pdd{z_i}{\zeta_j}{\zeta_k} 
              = \frac{1}{z_i-\zeta_k} \pd{z_i}{\zeta_j}
              + \frac{1}{z_i-\zeta_j} \pd{z_i}{\zeta_k} \\ 
   &\hskip7em
              - \frac{P''(z_i)}{P'(z_i)} \pd{z_i}{\zeta_j} 
              \pd{z_i}{\zeta_k} - \frac{1}{P'(z_i)} 
              \int_\beta^{z_i} \frac{P'(w)\,dw}{(w-\zeta_j)(w-\zeta_k)}
              \text{\quad for $j \ne k$}.
\end{align*}
\end{lemma}
\begin{proof}
We may assume that $z_i \ne \beta$ (else the results would be
trivially true.)  Proofs of parts 1--3 can be found in \cite[Lemmas
2.1 and 2.3]{Borcea-1}.  Part 4 can be established by writing
$P'(z)=(z-\zeta_j)Q(z)$ (with $Q$ independent of $\zeta_j$) and
calculating the second partial derivative (with respect to $\zeta_j$)
of
\begin{displaymath}
    0 = \int_\beta^{z_i} wQ(w)\,dw - \zeta_j \int_\beta^{z_i} Q(w)\,dw.
\end{displaymath}
Part 5 can likewise be established by writing
$P'(z)=(z-\zeta_j)(z-\zeta_k)Q(z)$ and calculating the mixed second
partial derivative.
\end{proof}

\section{Proof of Theorem~\ref{theorem2}} \label{section4} 

Recall that $P$ is a real polynomial of degree $n \ge 5$ that
satisfies properties \ref{A}--\ref{H}.  We will prove Theorem
\ref{theorem2} in two stages, first by estimating the roots and
critical points of $P$ with linear approximations, and then by
improving our estimates with quadratic approximations.

We begin by examining the conclusions that can be drawn from linear
approximations to the roots and critical points of $P$ with
\begin{proposition}\label{prop5}
If $I=(\Delta\beta, \Delta\zeta_1,\dots,\Delta\zeta_{n-1})$ is an 
improvement of $P$ and $\Delta t=||I||$, then
   $\Delta\beta = \Odt{2}$, 
   $\Delta\zeta_j = \Odt{2}$ for $j \le 2$,
   $\Re[\Delta\zeta_j] = \Odt{2}$ for $j \ge 3$ and
   $\sum_{j=3}^{n-1} \Delta\zeta_j = \Odt{2}$.
\end{proposition}
\begin{proof}
From the definition of $\Delta t$ we know that $\Delta\beta=\Odt{}$
and that each $\Delta\zeta_j=\Odt{}$.  Since $I$ is an improvement,
then each $|(\zeta_j-\beta)+(\Delta\zeta_j-\Delta\beta)| > r$, so
$|1+(\Delta\zeta_j-\Delta\beta)/(\zeta_j-\beta)| > r/|\zeta_j-\beta| =
1$, so using part 1 of Lemma \ref{lemma3} gives us that each
\begin{equation}\label{4.1}
      -\Re \left[ \frac{\Delta\zeta_j-\Delta\beta}{\zeta_j-\beta} \right]
      \le \Odt{2}.
\end{equation}
Recalling the expressions $E_k$ defined in \ref{2.1}, this gives us
two inequalities $E_1 \le \Odt{2}$ and $E_2 \le \Odt{2}$.  Adding
\eqref{4.1} for $j=3, \dots, n-1$ and recalling that $\zeta_j=a$ is
real for $j \ge 3$ gives us a third inequality $E_3 \le \Odt{2}$.

By part 1 of Lemma \ref{lemma4} the roots of our improved polynomial
\eqref{3.1} that originate from the simple roots $z_i$ of $P$ are
analytic functions of $\beta, \zeta_1, \dots, \zeta_{n-1}$, so each
such root is of the form $z_i+\Delta z_i$, with
\begin{displaymath}
   \Delta z_i = \pd{z_i}{\beta} \Delta\beta 
              + \sum_{j=1}^2 \pd{z_i}{\zeta_j} \Delta\zeta_j 
              + \pd{z_i}{\zeta_3} \sum_{j=3}^{n-1} \Delta\zeta_j
              +\Odt{2}.
\end{displaymath}
Note that each $\Delta z_i=\Odt{}$.  Since $I$ is an improvement, then
each $|z_i+\Delta z_i| \le 1$. If $|z_i|=1$, then $|1+\Delta z_i/z_i|
\le 1/|z_i| = 1$ so using part 1 of Lemma~\ref{lemma3} gives us $\Re[\Delta
z_i/z_i]\le \Odt{2}$ and thus we have inequalities $E_{i+3} \le \Odt{2}$
for $1=1, \dots, m$.

By property~\ref{G}, there are positive constants $c_k$ so that
$\sum_{k=1}^{m+3} c_k E_k=0$.  Since each $E_i\le \Odt{2}$, then each
\begin{displaymath}
   E_i = (-1/c_i) \sum_{i \ne k=1}^{m+3} c_k E_k \ge \Odt{2}
\end{displaymath}
so each $E_i=\Odt{2}$.

Thus we consider the system $\{E_k=\Odt{2}$ : $k=1, \dots, m+3\}$.  By
property~\ref{H}, the coefficient matrix of this system is of rank
$7$, and so solving this system shows that the values of our 7
``variables'' are all $\Odt{2}$.  Thus we can conclude that
$\Delta\beta=\Odt{2}$, that $\Delta\zeta_j=\Odt{2}$ for $j \le 2$, and
that $\sum_{j=3}^{n-1} \Delta\zeta_j=\Odt{2}$.

Suppose that $j \ge 3$ and note that $\sum_{k=3}^{n-1}
\Re[\Delta\zeta_k]=\Odt{2}$.  Since by property~\ref{E} we know that
$\zeta_j-\beta=a-\beta<0$, and since $\Delta\beta=\Odt{2}$, then
\eqref{4.1} implies that each $\Re[\Delta\zeta_j] \le \Odt{2}$, so
each
\begin{displaymath}
   \Re[\Delta\zeta_j] = -\sum_{j \ne k=3}^{n-1} \Re[\Delta\zeta_k]+\Odt{2} 
                      \ge \Odt{2}
\end{displaymath}
and so each $\Re[\Delta\zeta_j]=\Odt{2}$.  This finishes the proof of
Proposition \ref{prop5}.
\end{proof}

At this point, for $j \ge 3$ we know only that each
$\Delta\zeta_j=\Odt{}$.  We can improve this estimate by looking at
quadratic approximations to the roots and critical points of $P$ with
\begin{proposition}\label{prop6}
If $I=(\Delta\beta, \Delta\zeta_1,\dots,\Delta\zeta_{n-1})$ is an
improvement of $P$ and $\Delta t=||I||$, then each $\Delta\zeta_j =
\Odt{3/2}$.
\end{proposition}
\begin{proof}
Note that the hypotheses of Proposition \ref{prop5} are satisfied, so
we may use all of its conclusions.  In particular, we know that
$\Delta\zeta_j=\Odt{2}$ for $j \le 2$ and that
$\Re[\Delta\zeta_j]=\Odt{2}$ for $j \ge 3$, so to verify Proposition
\ref{prop6} we need only show that $\Im[\Delta\zeta_j]=\Odt{3/2}$ for
$j \ge 3$.  We will do this by repeating the calculations of
Proposition \ref{prop5}, but working now to $\Odt{3}$.

Since $I$ is an improvement, then each
$|1+(\Delta\zeta_j-\Delta\beta)/(\zeta_j-\beta)| > 1$, so using part 2
of Lemma \ref{lemma3} gives us
\begin{equation}
  \label{4.2}
      -\Re \left[ \frac{\Delta\zeta_j-\Delta\beta}{\zeta_j-\beta} \right]
      - \frac{1}{2} \left( \Im \left[ 
         \frac{\Delta\zeta_j-\Delta\beta}{\zeta_j-\beta} \right] \right) ^2 
      \le \Odt{3}.
\end{equation}
Define inequalities 1 and 2 by evaluating \eqref{4.2} for $j=1$ and
$j=2$ respectively.  By Proposition~\ref{prop5}, for $j \le 2$ we have
$\Delta\zeta_j-\Delta\beta=\Odt{2}$, so inequalities 1 and 2 are
$E_k+F_k \le \Odt{3}$ for $k=1, 2$.

Define inequality 3 to be the sum of \eqref{4.2} for $j=3, \dots,
{n-1}$.  Now each $\zeta_j-\beta=-(\beta-a)$ is real for $j \ge 3$,
and each $\Im[\Delta\zeta_j-\Delta\beta]=\Im[\Delta\zeta_j]$, so
inequality~3 can be written as $E_3+F_3 \le \Odt{3}$.

Each root of our improved polynomial is of the form $z_i+\Delta z_i$.
A quadratic approximation to $\Delta z_i$ includes terms of the form
$\Delta\beta\Delta\zeta_j=\Odt{3}$ and (for $j \le 2$ or $k \le 2$)
$\Delta\zeta_j\Delta\zeta_k=\Odt{3}$, which are absorbed into the
$\Odt{3}$ when we write
\begin{displaymath}
   \Delta z_i = \pd{z_i}{\beta} \Delta\beta 
            + \sum_{j=1}^2 \pd{z_i}{\zeta_j} \Delta\zeta_j
            + \pd{z_i}{\zeta_3}\sum_{j=3}^{n-1} \Delta\zeta_j
            +\frac{1}{2} \sum_{j=3}^{n-1} \sum_{k=3}^{n-1} 
                 \pdd{z_i}{\zeta_j}{\zeta_k} \Delta\zeta_j \Delta\zeta_k
            +  \Odt{3}.
\end{displaymath}
Note that Proposition \ref{prop5} implies that each $\Delta
z_i=\Odt{2}$.  Now $\zeta_j=a$ for $j \ge 3$, so Lemma \ref{lemma4}
shows that for $j \ge 3$ and $k \ge 3$ we have
\begin{displaymath}
\pdd{z_i}{\zeta_j}{\zeta_k} - \pdd{z_i}{\zeta_3}{\zeta_4} = 
   \begin{cases} \displaystyle 
       \frac{1}{P'(z_i)} \int_\beta^{z_i} \frac{P'(w)\,dw}{(w-a)^2}
           &\text{if $j=k$}\medskip\\
       0   &\text{if $j \ne k$.}
   \end{cases}
\end{displaymath}
From Proposition \ref{prop5} we know that
\begin{displaymath}
   \sum_{j=3}^{n-1} \sum_{k=3}^{n-1} \Delta\zeta_j \Delta\zeta_k
      = \left( \sum_{j=3}^{n-1} \Delta\zeta_j \right)^2 = \Odt{4}
\end{displaymath}
and that $(\Delta \zeta_j)^2 = -(\Im[\Delta\zeta_j])^2 + \Odt{3}$ for $j
\ge 3$, so
\begin{align*}
   \sum_{j=3}^{n-1} \sum_{k=3}^{n-1} \pdd{z_i}{\zeta_j}{\zeta_k}
         \Delta\zeta_j \Delta\zeta_k
   &=  \sum_{j=3}^{n-1} \sum_{k=3}^{n-1} 
       \left( \pdd{z_i}{\zeta_j}{\zeta_k} - \pdd{z_i}{\zeta_3}{\zeta_4} \right)
         \Delta\zeta_j \Delta\zeta_k +\Odt{4} \\
   &=  \frac{-1}{P'(z_i)} \int_\beta^{z_i} \frac{P'(w)\,dw}{(w-a)^2}
       \sum_{j=3}^{n-1} (\Im[\Delta\zeta_j])^2 +  \Odt{3}.
\end{align*}

Recall that each $|z_i+\Delta z_i| \le 1$ and that each $\Delta
z_i=\Odt{2}$. If $|z_i|=1$, then $|1+\Delta z_i/z_i| \le 1$, so using
part~1 of Lemma~\ref{lemma3} gives us $\Re[\Delta z_i/z_i] \le
\Odt{4}$ and so $E_{i+3}+F_{i+3} \le \Odt{3}$ for $i=1,\dots,m$.

Thus we have $E_k+F_k \le \Odt{3}$ for $k=1,\dots,m+3$.  From
property~\ref{G} we have positive constants $c_k$ so that
$\sum_{j=3}^{n-1} (\Im[\Delta\zeta_j])^2 = \sum_{k=1}^{m+3}
c_k(E_k+F_k) \le \Odt{3}$ and thus $\Im[\Delta\zeta_j] = \Odt{3/2}$
for $j \ge 3$, which completes the proof of Proposition \ref{prop6}.
\end{proof}

Using our estimates from Propositions \ref{prop5} and \ref{prop6}, we
can now write the 
\begin{proof}[Proof of Theorem \ref{theorem2}]
Let $I=(\Delta\beta, \Delta\zeta_1,\dots,\Delta\zeta_{n-1})$ be any
improvement of $P$, and let $\Delta t=||I||$.  From Propositions
\ref{prop5} and \ref{prop6}, we know that $\Delta\beta=\Odt{2}$ and
that each $\Delta\zeta_j=\Odt{3/2}$, so
\begin{displaymath}
   \Delta t = ||I|| 
    = \left( |\Delta\beta|^2 + \sum_{j=1}^{n-1} |\Delta\zeta_j|^2 \right)^{1/2}
     = \Odt{3/2}.
\end{displaymath}
Thus there is a constant $K>0$ so that $\Delta t \le K(\Delta t)^{3/2}$, so
$||I||=\Delta t \ge 1/K^2$.
\end{proof}

\section{Calculations}\label{section5} 

Recall that Theorem~\ref{theorem2} implies that any real polynomial of
degree $n \ge 5$ that satisfies Properties \ref{A}--\ref{H} is locally
extremal.  In this section, we show how to find locally extremal
polynomials by describing how to construct polynomials that satisfy
our properties~\ref{A}--\ref{H}.  Maple code for these calculations
can be found in the source files for this paper, located at
http://arxiv.org/abs/math/0505424 [v3].

Recall that $P$ is to be a real polynomial of degree $n \ge 5$.  By
property~\ref{E} we must have $P'(z)=(z-a)^{n-3}(z^2+bz+c)$ for some
real numbers $a$, $b$ and $c$.  By property~\ref{C}, we know that
$P(z)=\int_\beta^z P'(w)\,dw$, so the coefficients of $P$ are
polynomials in $\{\beta, a, b, c\}$.

Applying property~\ref{D} to the critical points of $P$ (the real
number $a$ and the complex roots of $z^2+bz+c$), we get our first
equation $\beta^2+b \beta+c=(\beta-a)^2$.

Note that property~\ref{G} implies that the rows of the linear system
$\{E_k =0: k=1,\dots,m+3\}$ are linearly dependent.  Since
property~\ref{H} states that the coefficient matrix of this system is
of rank 7, this means that there must be at least 8 equations in the
system, so $m \ge 5$ and thus we must look for polynomials with at
least 5~roots on the unit circle.

If $n$ is odd, we will look for polynomials with three pairs of
complex conjugate roots on the unit circle.  Each such pair will be
the roots of a quadratic of the form $z^2+d_i z+1$, so the remainders
upon dividing $P$ by each of these three quadratics are linear
polynomials with both coefficients equal to $0$.  This generates an
additional six equations in the seven variables $\{\beta, a, b, c,
d_1, d_2, d_3\}$. Thus we have a nonlinear system of 7 equations in 7
unknowns.

If $n$ is even, we will look for polynomials with two pairs of complex
conjugate roots on the unit circle and a root at $-1$.  The conjugate
roots generate an additional 4 equations in the six variables
$\{\beta, a, b, c, d_1, d_2\}$ (as above), and the equation
$P(-1)=0$ generates a sixth equation. Thus we have a nonlinear system
of 6 equations in 6 unknowns.

Thus in either case we get a nonlinear system of equations, with the
same number of equations as unknowns, so we can try to solve this
system.  Note that there may be more than one solution, so we will
need to choose the ``correct'' one.  For each such solution, we verify
properties~\ref{A}--\ref{H} by numerically checking the following
assertions.  (Details of these computations may be found in the Maple
code referenced above.)

Property~\ref{A}: The maximum modulus of the roots of $P$ is equal to
1 (to the accuracy calculated.)

Property~\ref{B}: The minimum distance between any two roots of $P$ is
greater than~$0.1$.  (This shows that all roots of $P$ are simple.)

Property~\ref{C}: We have $0.7<\beta<0.9$.  (We know that
$P(\beta)=0$, since $P(z)=\int_\beta^z P'(w)\,dw$ by construction.)

Property~\ref{D}: The distances between $\beta$ and the critical
points of $P$ are all equal (to the accuracy calculated), and this
common distance is greater than $0.9$.

Property~\ref{E}: The quantity $\beta-a>0.9$.  (Note that the critical
point $a$ is real and of order $n-3$ by construction.)

Property~\ref{F}: If we define 
\begin{displaymath}
R=\max \left\{\min_{1 \le j \le n-1} |z_i-\zeta_j| : z_i \ne \beta \right\}
        \text{ and }
r=\min_{1 \le j \le n-1} |\beta-\zeta_j|, 
\end{displaymath}
then $r<0.97$ and $r-R>0.02$.

Property~\ref{G}: The linear system given by the sums has a solution
in which every $c_k>0.3$.

Property~\ref{H}: The seventh largest singular value of the
coefficient matrix is greater than $0.04$.

Once we have verified properties~\ref{A}--\ref{H} for a specific
polynomial $P$, we know by Theorem \ref{theorem2} that $P$ is locally
extremal and we are done.

\section{Proof of Theorem \ref{theorem1}}\label{section6} 

For the values of $\beta$ and $P'(z)$ given below, one can verify that
the polynomials $P(z)=\int_\beta^z P'(w)\,dw$ satisfy
Properties~\ref{A}--\ref{H}, and thus by Theorem~\ref{theorem2} are
locally extremal.  (Details of these computations may be found with
the Maple code referenced above.)

For $n=8$, we take $\beta=0.7290857513$ and
\begin{displaymath}
P'(z) = (z+0.2035409790)^5\, (z^2 - 0.5410836525 z + 0.7327229666).
\end{displaymath}

For $n=9$, we take $\beta=0.7145672829$ and
\begin{displaymath}
P'(z) = (z+0.2157115753)^6\, (z^2 - 0.8021671918 z + 0.9280147829).
\end{displaymath}

For $n=12$, we take $\beta=0.8403619619$ and
\begin{displaymath}
P'(z) = (z+0.1155828545)^9\, (z^2 - 0.4090272613 z + 0.5513532168).
\end{displaymath}

For $n=13$, we take $\beta=0.8275325585$ and
\begin{displaymath}
P'(z) = (z+0.1246203379)^{10}\, (z^2 - 0.5415308686 z + 0.6699194279).
\end{displaymath}

For $n=14$, we take $\beta=0.8158105092$ and
\begin{displaymath}
P'(z) = (z+0.1304708647)^{11}\, (z^2 - 0.6885970233 z + 0.7916663399).
\end{displaymath}

For $n=15$, we take $\beta=0.7999767588$ and
\begin{displaymath}
P'(z) = (z+0.1400336168)^{12}\, (z^2 - 0.8389864647 z + 0.9148263642).
\end{displaymath}

For $n=19$, we take $\beta=0.8684432238$ and
\begin{displaymath}
P'(z) = (z+0.0923361850)^{16}\, (z^2 - 0.6503807257 z + 0.7337221736).
\end{displaymath}

For $n=20$, we take $\beta=0.8570396874$ and
\begin{displaymath}
P'(z) = (z +0.0982636528)^{17}\, (z^2 - 0.7563752823 z + 0.8263310816).
\end{displaymath}

For $n=26$, we take $\beta=0.8817716692$ and
\begin{displaymath}
P'(z) = (z +0.0797127446)^{23}\, (z^2 - 0.7969496845 z + 0.8496586550).
\end{displaymath}

These verifications complete the proof of Theorem~\ref{theorem1}.


\end{document}